\documentclass[12pt, reqno]{amsart}

\usepackage{amsmath}
\usepackage{amsfonts}
\usepackage{amssymb}

\begin{document}
\title[On the asymptotics of intergals]{On the asymptotics of intergals arising \\ in the study of the wave motion in lattices}
\author{Nadezhda I.~Aleksandrova}
\address{Chinakal Institute of Mining, Siberian Branch, Russian Academy of Sciences, Novosibirsk, Russia}
\email{nialex@misd.ru}

\begin{abstract}
We present a method for finding the asymptotics of integrals arising in solid mechanics.
\end{abstract}
\maketitle
\thispagestyle{empty}
\section{Dynamics of a block environment}

Until recently, in geomechanics and geophysics, the theory of deformations of a rock mass as a homogeneous medium was widely used.  
The dynamics of such a homogeneous medium is described by a well-developed linear theory of elastic wave propagation.  
For the first time, a serious reason for revising the prevailing views was given in [1]. 
There, the need was indicated to take into account the block structure of rocks 
in mathematical models intended for geomechanics and seismics.
The main idea proposed in [1] is to consider
the rock mass as a system of blocks of different scale levels nested into each other, 
connected by layers consisting of weaker, fractured rocks.
 
In the simplest case, the dynamics of a block medium is studied in the pendulum approximation, 
when it is believed that the blocks are incompressible, and all deformations and displacements 
occur due to the compressibility of the interlayers only.
In this case, the block medium can be modeled as a lattice of point-masses connected to each other by springs.
Within the framework of this model, we study the antiplane deformation of a two-dimensional square lattice consisting 
of point-masses, having the same masse $M$, connected by springs of length $L$, having the same stiffness $k$.
The values $M$, $L$, and $k$ are taken as units.

\section{Equations of motion and their solutions}

The equations of motion of point-masses have the form
$$
\ddot{u}_{m,n}  =u_{m+1,n} +u_{m-1,n} +u_{m,n+1} +u_{m,n-1} -4u_{m,n} +Q(t),
$$
where $u_{m,n}$ is the displacement of the mass in the direction orthogonal to the lattice plane; 
$m,n$ are the mass numbers in the directions of the axes $x,y$; 
$Q(t)=Q_0\sin(\omega _*t)H(t)\delta(n)\delta(m)$ is the sinusoidal load of amplitude $Q_0$ and frequency $\omega_*$
applied to the point with coordinates $(0,0)$; 
$H$ is the Heaviside function;
$\delta$ is the Dirac function.

Applying the Laplace transform with respect to time $t$ and the discrete Fourier transform in the variables $m,n$ 
to the above equations, we obtain the Laplace--Fourier transform of the solution:
$$
u^{LF_mF_n}=\frac{Q^{LF_mF_n}}{p^2+2(2-\cos q_x-\cos q_y)},\quad Q^{LF_mF_n}=\frac{Q_0\omega_*}{p^2+\omega^2_*}.
$$
Here $L$ denotes the Laplace transform with parameter $p$; 
$F_m$, $F_n$ denotes the discrete Fourier transforms in $m$, $n$ with parameters $q_x$, $q_y$.

Inverting thus obtained expression for $u^{LF_mF_n}$ with respect to $n$, we get
$$
u^{LF_m}_n =\frac{Q_0\omega_*(B-\sqrt{B^2-1})^{|n|}}{2(p^2+\omega^2_*)\sqrt{B^2-1}},
\qquad\mbox{where}\quad
B =\frac{p^2}{2}+2-\cos q_x.
$$
 
For $u^{LF_m}_n$, it is not possible to find  explicitly the inverse Laplace and Fourier transforms.
For this reason, we look for the asymptotic behavior of $u_{m,n}$ for $t\to\infty$.

The frequency $\omega=2$ is resonant for the lattice of point-masses under consideration [2]. 
Put by definition $\omega_*=2$ and $p=s+2i$.
Then
$$
u^{LF_m}_n\sim \frac{Q_0(-1)^{n+1}\exp{(i|q_x n|)}}{4s\sqrt{\sin^2 q_x+4is\cos q_x}}
\quad\mbox{as}\ s\to 0.
$$ 
Here and below, the formula $v(z)\sim w(z)$ for $z\to z_0$ means that $\lim_{z\to z_0} [v(z)-w(z)]=0$.

Using the inversion formula for the Fourier transform, we obtain
$$
u^{L}_{m,n}\sim\frac{Q_0(-1)^{n+1}}{4\pi s}\int_{0}^{\pi}
\frac{ \cos(q_x m)\exp (iq_x |n|)}{\sqrt{\sin^2 q_x + 4is \cos q_x}}\, dq_x\quad\mbox{as}\ s\to 0. \eqno(1)
$$
 
The asymptotics of integrals (1) as $s\to 0$, which corresponds to $t\to \infty $ in the space of originals, has the following form:

a) if $m+n$ is even, then
$$
u^L_{0,0}(s)\sim -\frac{Q_0}{4\pi s}\ln \frac{4}{s},\quad u_{0,0}(t)\sim \frac{Q_0}{4\pi}[\ln(4t)+\gamma], 
$$
$$
u^L_{m,m}(s)\sim\frac{Q_0(-1)^m}{4\pi s}[\ln(s|m|)+\gamma ],\quad m\neq 0, \eqno(2)
$$
$$
u_{m,m}(t)\sim \frac{Q_0(-1)^{m+1}}{4\pi}\ln \frac{t}{|m|}, \quad m\neq 0, 
$$
$$
u^L_{m,n}(s)\sim\frac{Q_0(-1)^n}{4\pi s}[\ln(s|m^2-n^2|)+2\gamma],\quad |m|\neq |n|, \eqno(3)
$$
$$   
u_{m,n}(t)\sim\frac{Q_0(-1)^{n+1}}{4\pi }[\ln \frac {t}{|m^2-n^2|}-\gamma],\quad |m|\neq |n|, 
$$
where $\gamma\approx 0.577$ is the Euler constant;

b) if $m+n$ is odd, then
$$
u^L_{m,n}(s)\sim\frac{iQ_0(-1)^{n+1}}{16s},\quad u_{m,n}(t)\sim\frac{iQ_0(-1)^{n+1}}{16}. \eqno(4)
$$

We emphasize that formulas (2)--(4) are more precise than similar formulas found in [2].

\section{Finding the asymptotics of an integral}

We demonstrate the method of finding asymptotics of integrals (1) by the example of 
finding asymptotics of a simpler integral
$$
\int_{0}^{\pi}\frac{dq}{\sqrt{\sin^2 q+is}}.
$$
This integral is similar to (1) for $m=n=0$. 
The problem of finding asymptotics of this integral was published by the author in [3].
A solution to this problem, which is different from the solution discussed below, is also published in
``The American Mathematical Monthly,'' see [3].

Put by definition
$$
h(q, s)=\frac{1}{\sqrt{\sin^2 q+is}}.
$$
Since 
$$
h\biggl(\dfrac{\pi}{2}+q,s\biggr)=h\biggl(\dfrac{\pi}{2}-q,s\biggr)
$$ 
for all $s>0$ and $q\in [0,\pi]$, then
$$\int_0^\pi\frac{dq}{\sqrt{\sin^2 q+is}}=2\int_0^{\pi/2} h(q,s)\,dq
= 2\int_0^\varepsilon h(q,s)\,dq+2\int_\varepsilon^{\pi/2} h(q,s)\,dq. \qquad { } \eqno(5)
$$
Here $\varepsilon>0$ is a fixed sufficiently small number.
Using the notation
$h_1(q,s)=\text{Re\,} h(q,s)$ and $h_2(q,s)=\text{Im\,} h(q,s)$,
we get $h(q,s)=h_1(q,s)+ih_2(q,s),$
$$
h_1(q,s) =\sqrt{
\frac{\sqrt{\sin^4 q+s^2}+\sin^2 q}{2(\sin^4 q+s^2)}
},\qquad
h_2(q,s) =\sqrt{
\frac{\sqrt{\sin^4 q+s^2}-\sin^2 q}{2(\sin^4 q+s^2)}
}.
$$

We find separately the asymptotics of the real and imaginary parts of the both integrals on the right-hand side of formula (5):
$$
 2\int_\varepsilon^{\pi/2} h_1(q,s)\,dq\underset{s\to+0}{\sim}
2\int_\varepsilon^{\pi/2} \frac{dq}{\sin q}
 =2\ln \biggl(\frac{1}{\tan (\varepsilon/2)}\biggr)
\underset{\varepsilon\to+0}{\sim}\ln\biggl(\frac{4}{\varepsilon^2}\biggr),\qquad { } \eqno(6)
$$
$$
2\int_0^\varepsilon h_1(q,s)\,dq=
2\int_0^\varepsilon \sqrt{
\frac{\sqrt{\sin^4 q+s^2}+\sin^2 q}{2(\sin^4 q+s^2)}
}\, dq 
$$
$$
\underset{\varepsilon\to+0}{\sim}
2\int_0^\varepsilon \sqrt{
\frac{\sqrt{q^4 +s^2}+ q^2}{2(q^4 +s^2)}
}\,dq\overset{\boxed{1}}
= 2\int_0^{\varepsilon/\sqrt{s}}\sqrt{
\frac{\sqrt{y^4+1}+y^2}{2(y^4+1)}
}\,dy
$$
$$
\overset{\boxed{2}}{=}\int_1^{\varphi(\varepsilon/\sqrt{s})}
\frac{dz}{\sqrt{z^2-1}}=
\ln\biggl(\sqrt{\varphi^2\biggl(\frac{\varepsilon}{\sqrt{s}}\biggr)-1}
+\varphi\biggl(\frac{\varepsilon}{\sqrt{s}}\biggr)
\biggr)
\underset{s\to+0}{\sim}
\ln\biggl(4\frac{\varepsilon^2}{s}
\biggr)\qquad { } \eqno(7)
$$
$\bigl($in (7), the following substitutions are used:
$q=y\sqrt{s}$ in order to get equality $\boxed{1}$
and $\sqrt{y^4+1}+y^2=z$ or, for short, $z=\varphi(y)$ in order to get equality $\boxed{2}\bigr)$;
$$
2\int_\varepsilon^{\pi/2} h_2(q,s)\,dq\underset{s\to+0}{\sim}
\int_\varepsilon^{\pi/2} \frac{s\, dq}{\sin^3 q}\underset{s\to+0}{\to}0,\eqno(8)
$$
$$
 2\int_0^\varepsilon h_2(q,s)\,dq=
2\int_0^\varepsilon \sqrt{
\frac{\sqrt{\sin^4 q+s^2}-\sin^2 q}{2(\sin^4 q+s^2)}
}\, dq 
$$
$$
\underset{\varepsilon\to+0}{\sim}
2\int_0^\varepsilon \sqrt{
\frac{\sqrt{q^4 +s^2}- q^2}{2(q^4 +s^2)}
}\,dq\overset{\boxed{3}}{=}2\int_0^{\varepsilon/\sqrt{s}}\sqrt{
\frac{\sqrt{y^4+1}-y^2}{2(y^4+1)}
}\,dy
$$
$$
\overset{\boxed{4}}{=} -\int_1^{\psi(\varepsilon/\sqrt{s})}
\frac{dz}{\sqrt{1-z^2}}=
-\arcsin z\biggl\vert_1^{\psi(\varepsilon/\sqrt{s})}
\underset{s\to+0}{\sim}
\frac{\pi}{2}\eqno(9)
$$
$\bigl($in (9), the following substitutions are used:
$q=y\sqrt{s}$ in order to get equality $\boxed{3}$
and $\sqrt{y^4+1}-y^2=z$ or, for short, $z=\psi(y)$ in order to get equality $\boxed{4}\bigr)$.

Substituting formulas (6)--(9) into (5), we obtain
$$
{}\qquad\qquad
\int_0^\pi\frac{dq}{\sqrt{\sin^2 q+is}}\sim \ln \biggl(\frac{16}{s}\biggr)+i\frac{\pi}{2}\quad\text{as}\quad s\to +0. 
$$

The above asymptotics (2)--(4) of integral (1) are obtained in the same way.

\end{document}